\documentclass[12pt,lot, lof]{amsart}
\author[Mohammad F. Tehrani]{Mohammad Farajzadeh Tehrani}
\address{Simons center for geometry and physics}
\email{mtehrani@scgp.stonybrook.edu}

\title[Automorphism group of Batyrev Calabi-Yau threefolds]
{Automorphism group of Batyrev Calabi-Yau threefolds}

\usepackage{epsfig}
\usepackage{amsfonts,amssymb,amsmath}   
\usepackage{amsthm}
\usepackage{latexsym}
\usepackage{graphicx}
\usepackage[all]{xy}
\usepackage{verbatim}           
\usepackage{multirow}
\usepackage{longtable}
\usepackage{booktabs}
\usepackage{hyperref}
\hypersetup{bookmarksnumbered}
\hypersetup{colorlinks,bookmarksnumbered}
\usepackage{enumitem}

\numberwithin{equation}{section}

\newtheorem{theorem}{Theorem}[section]
\newtheorem{lemma}[theorem]{Lemma}

\newtheorem{proposition}[theorem]{Proposition}
\newtheorem{conjecture}[theorem]{Conjecture}

\theoremstyle{definition}
\newtheorem{definition}[theorem]{Definition}
\newtheorem{remark}[theorem]{Remark}
\newtheorem{example}[theorem]{Example}

\def\mc{\mathcal}

\def\ep{\epsilon}
\def\ra{\rightarrow}
\def\De{\Delta}

\def\om{\omega}

\def\si{\sigma}
\def\Si{\Sigma}

\def\ov#1{\overline{#1}}
\def\tn#1{\textnormal{#1}}

\def\bEqu#1{\begin{equation}\label{#1}}
\def\eEqu{\end{equation}}
\def\bsEq{\begin{equation*}}
\def\esEq{\end{equation*}}
\def\bDef#1{\begin{definition}\label{#1}}
\def\eDef{\end{definition}}
\def\bThm#1{\begin{theorem}\label{#1}}
\def\eThm{\end{theorem}}
\def\bLem#1{\begin{lemma}\label{#1}}
\def\eLem{\end{lemma}}
\def\bRem#1{\begin{remark}\label{#1}}
\def\eRem{\end{remark}}
\def\bExa#1{\begin{example}\label{#1}}
\def\eExa{\end{example}}
\def\bPro#1{\begin{proposition}\label{#1}}
\def\ePro{\end{proposition}}
\def\bFig#1{\begin{figure}\label{#1}}
\def\eFig{\end{figure}}
\def\bProof{\begin{proof}}
\def\eProof{\end{proof}}

\def\sing{\tn{sing}}

\def\co{\mc{O}}

\def\ck{\mc{K}}

\def\ck{\mc{K}}

\def\R{\mathbb R}
\def\C{\mathbb C}
\def\Z{\mathbb Z}
\def\Q{\mathbb Q}
\def\P{\mathbb P}

\begin{document}
\begin{abstract}
In this paper, we will prove that all Batyrev Calabi-Yau threefolds, arising from a small resolution of a generic hyperplane section of a reflexive Fano-Gorenstein fourfold, have finite automorphism group. 
Together with Morrison conjecture, this suggests that Batyrev Calabi-Yau threefolds should have a polyhedral K\"{a}hler (ample) cone.
\end{abstract}
\maketitle
\section{Introduction}\label{sec:intro}

The construction of Calabi-Yau threefolds, that brought by far the largest amount of examples is the construction of Batyrev \cite{B1} and its generalization \cite{BB}.
Starting from a four dimensional reflexive polytop $\De$, let $X_{\De}$ be the corresponding toric 4-fold. Then a generic section of anti-canonical line bundle is a singular Calabi-Yau threefold $Z_\De$. 
There is always a crepant resolution $\ov{Z}_\De$ of $Z_\De$, although it is not unique. 
This resolution is induced by a maximal projective crepant partial (MPCP) desingularization $\ov{X}_{\De}$ of $X_{\De}$. 
This way, we obtain a smooth Calabi-Yau 3-fold, which we call a \textsf{Batyrev Calabi-Yau 3-fold}. 
Reflexive polytops in dimension four are classified and the number is 473,800,776. 
Among them, there exist at least 30000 different topological types of Calabi-Yau threefolds, i.e. with different Hodge numbers. 
In the same paper, Batyrev examined the mirror symmetry conjecture on this set of manifolds and proved that the dual Calabi-Yau families, obtained from the dual polytops, $\De$ and $\De^*$, have mirror hodge numbers. This construction, generalizes the previously proposed mirror construction of quintic threefold,  due to physicists \cite{CDGP}.\\
 
Given a smooth Calabi-Yau 3-fold $Z$ \footnote{We define a smooth C.Y 3-fold to be a simply-connected smooth projective 3-fold with trivial canonical bundle.}, 
its K\"{a}hler cone $\ck_Z \subset H^{1,1}(Z,\R)$, is the set of cohomology classes which can be represented by a K\"{a}hler form. 
Under the isomorphism 
$$H^2(Z,\R)=H^{1,1}(Z,\R)=\tn{Pic}(Z)\otimes_\Z \R\cong  N^1(Z)\otimes_\Z \R,$$
$\ck_Z$ is isomorphic to the ample-cone; and its closure, nef cone, is dual to the closure of cone of effective curves, $\ov{\tn{NE}}(Z,\R)\subset N^1(Z)\otimes_\Z \R \cong H_2(Z,\R)$. 
We will denote the compactified K\"{a}hler cone by $\ov\ck_Z$. 
The group of automorphims of $Z$, $\tn{Auto}(Z)$, acts on $\ck_Z$.

The K\"{a}hler cone of $Z$ can have a rather complicated geometry, see \cite{W}. Away from cohomology classes with triple intersection zero, $\ov\ck_Z$ is locally rational polyhedral, but the rational external rays may accumulate toward these points. 
Therefore, understanding the geometry of K\"{a}hler cone, and determining whether it is generated by a finite set of rational rays is a hard question.
The cone conjecture of \cite{M2} predicts that while $\ov\ck_Z$ may have infinitely many external rays, there will only be finitely many edges up to the action of $\tn{Auto}(Z)$.

\begin{conjecture}[The Cone Conjecture]
Let $Z$ be a Calabi-Yau manifold, and suppose that $h^{2,0}(Z)=0$. Let $\Lambda = H^2(Z,\Z)/\tn{torsion}$, and let $\ov\ck_Z^\Q$ be the convex hull of $\ov\ck_Z\cap \Lambda\otimes_\Z \Q$. Then there exists a rational polyhedral cone $\Pi\subset \ov\ck_Z^\Q$ such that $\tn{Auto}(Z)\cdot \Pi=\ov\ck_Z^\Q$.
\end{conjecture}

A non-trivial case of the cone conjecture, over Calabi-Yau threefolds which are fiber products of generic rational elliptic surfaces with section, has been checked by Grassi and Morrison \cite{GM}. 
Borcea \cite{Bor} has checked this conjecture for desingularized Horrocks-Mumford quintics. 
For smooth Calabi-Yau hypersurfaces in smooth Fano fourfolds, Kollar \cite[Appendix B ]{BK} proved that the K\"{a}hler cone is rational polyhedral. 
The K3 version of this conjecture is proved by Kovacs \cite{Kov}. 

 Let $\mc{D}(Z)=\Lambda_\R+i\ck$ be the associated tube domain. 
For every $Z$ for which the cone conjecture holds, the K\"{a}hler moduli space $\mc{D}(Z)/\Lambda \rtimes \tn{Auto}(Z)$, will admit both a Satake-Baily-Borel (toroidal) compactification, and also a smooth Mumford compactification; see \cite{M2}. 
Thus, the cone conjecture plays a vital role in the $A$-side of mirror symmetry. 

Let $Z$ be a Calabi-Yau threefold as before and let $c_2(TZ) \in H^4(Z,\Z) \cong N_1(Z)$ be the second Chern class of tangent bundle. 
By \cite[Theorem 1.1]{Mi}, for every $\om \in \ck_Z$, $\om\cdot c_2(TZ) >0$; therefore $c_2(TZ)$, as a linear function on $H^2(Z,\R)$, is non-negative on $\ov\ck_Z$. 
If $c_2(TZ)$ is strictly positive on the closure of K\"{a}hler cone of a Calabi-Yau threefold $Z$,  then $\tn{Auto}(Z)$ is finite; see \cite{W1} or the argument before the proof of Theorem~\ref{thm:main}.
Following theorem is the main result of this paper.
\bThm{thm:main}
Let $\ov{Z}_\De$ be a generic Batyrev Calabi-Yau threefold, then $c_2(T\ov{Z}_\De)$, restricted to $\ov\ck_{\ov{Z}_\De}$, is strictly positive.
\eThm

Therefore, we conclude that all Batyrev Calabi-Yau threefolds have finite automorphism group; together with the cone conjecture, this implies that the K\"{a}hler cone of all Batyrev Calabi-Yau threefolds has to be rational polyhedral. Note that if $X_{\De}$ is smooth, by the aforementioned result of Kollar, $\ov\ck_{\ov{Z}_\De}$ is rational polyhedral. We expect that some similar techniques should be applicable to the general case.

\bRem{rem:rigidity}
In \cite{W}, Wilson shows that the (geometry) of K\"{a}ler cone is almost constant in families. 
In fact, he proves \cite[Main theorem]{W} that outside countably many complex codimension one submanifolds (those Calabi-Yau threefolds containing a smooth elliptic ruled surface), 
K\"{a}hler cone will be invariant under the deformation of complex structure. 
Therefore, if it is rational polyhedral somewhere in a chamber, it would be rational polyhedral all over that chamber. 
This means that the genericity assumption of Theorem~\ref{thm:main} is not truly necessary.
\eRem

\subsection{Acknowledgment}
I would like to thank Gang Tian and P.M.H.Wilson for many helpful discussions.

\section{Review of Batyrev construction}\label{sec:basic-facts}
Let $M \cong \Z^n$ be a free abelian group of rank $n$; $M$ is a lattice in vector space $M_{\R} = M\otimes \R$ and $N = Hom(M,\Z)$ is the dual lattice in $N_{\R}=N\otimes \R$. 
Let $\De$ be an integral polytop, i.e. convex hull of finitely many lattice points, in $M$ including $0$ as an interior point. 
Every such polytop corresponds to a possibly singular toric algebraic variety $X_{\De}$ and an ample Cartier divisor $D_{\De}$ on it. $D_{\De}$ is the sum of toric divisors corresponding to the codimension one faces of $\De$.

\bDef{def:dual}
For every arbitrary rational polytop $\Delta$ in $M_{\Q}$, i.e. with vertices on $M\otimes_\Z \Q$, containing zero, its \textsf{dual} is defined by
$$ \Delta^*=\left\{y\in N_\Q \mid  \left\langle x,y\right\rangle 
\geq -1, \tn{for all}~ x\in \De \right\} \subset N_\Q,$$
where $\left\langle ,\right\rangle: M\times N \to \Z$ is the natural non-degenerate pairing between $N$ and $M$.
\eDef

\bDef{def:distance}
Let $H$ be a rational affine hyperplane in $M_\Q$ and $p\in M$ be an arbitrary integral point. 
Assume $H$ is affinely generated by integral points $H\cap M$, i.e. there exists $n_{H}\in N$ and an integer $c$ such that
$$ H=\left\{x \in M_\Q \mid \left\langle x,n_H\right\rangle=c \right\}.$$
Then, $\left|c-\left\langle p,n_H\right\rangle\right|$ is defined to be the \textsf{integral distance} between $H$ and $p$.
\eDef

\bDef{def:reflexive}
Let $0\in \De \subset M$ be a integral polytop. $(\De,M)$ is said to be a  \textsf{reflexive pair}, 
if the integral distances between $0\in M$ and all codimension-one faces of $\De$ are equal to 1. 
In this case, $(\De^*,N)$ is a reflexive pair as well.
\eDef

Through the rest of this paper, we set $(\De,\De^*)$ to be a pair of 4-dimensional reflexive polytops. 
Then, $X_{\De}$ corresponds to the fan given by the cones over faces of $\De^*$; and vice versa. 
We will denote the fan of $X_{\De}$ by $\Sigma$ and the fan of $X_{\De^*}$ by $\Sigma^*$. 
For a reflexive polytop $\De$, $D_{\De}$ is the anti-canonical divisor.

\bPro{pro:main}\textnormal(\cite[Theorem 4.1.9]{B1}\textnormal)
Given a 4-dimensional reflexive polytop $\De$,
\begin{itemize}
\item $D_\De$ is the section of the very ample line bundle $\co_\De(1)$ and $X_\De$ is a Fano-Gorenstein variety with canonical singularities;
\item the linear system $\left|-K_{X_\De}\right|$ is base-point free and its generic section $Z_\De$ is a 
Calabi-Yau model (it has a trivial canonical bundle) with at most canonical singularities;
\item $X_\De$ is smooth in codimension-one and its singular locus is supported on the collection of toric subvarieties $V_\Theta$, corresponding to some faces $\Theta$ of $\De$, of dimension 0,1,2 (point, curve and surface singularities, respectively);
\item $Z_\De$ is also smooth at codimension-one. $Z_\De^\sing$ has curve and point components corresponding to intersection points with the singular locus of $X_\De$; 
\item and finally, with $\Theta$ as before, if $\dim \Theta=0$, generic $Z_{\De}$ does not pass through the corresponding singular point; if $\dim \Theta =1$, $Z_\De \cap V_{\Theta}$ is a set of $d(\Theta):=\tn{length of } \Theta$ singular points on $Z_{\De}$; and if $\dim \Theta=2$, $Z_\De\cap V_\Theta$ is an irreducible curve of singularities on $Z_\De$. 
\end{itemize}
\ePro

Given a fan $\Si$ in $N$, a \textsf{refinement} $\Si'$ of $\Si$, given by subdividing the cones of $\Si$ into smaller ones, corresponds to a morphism of toric varieties $X_\Si' \to X_\Si$. 
If $\Si$ is given by cones over the faces of a polytop $\De^*$, a triangulation of the faces of $\De^*$ , i.e. a subdivision of faces of $\De^*$ into smaller simplecies, gives a refinement of $\Si$; but not every such triangulation results in a projective variety. 
We call a triangulation of the faces of $\De^*$ \textsf{projective}, if the corresponding variety $X_{\Si'}$ is projective, and we call it \textsf{maximal}, if it can not be refined further more. 

\bPro{pro:refinements}\textnormal(\cite{B1}\textnormal)
Given a 4-dimensional reflexive polytop $\De$, 
there exists a maximal projective triangulation of $\Si$. 
Moreover, every such triangulation corresponds to a maximal projective crepant partial desingularization (MPCP) $\pi\colon \ov{X}_\De=X_{\Si'} \to X_\De$ such that
\begin{itemize}
\item $\ov{X}_\De$ is smooth in codimension three and has only $\Q$-factorial singular points. In fact $\Si'$ will be simplicial and $\ov{X}_\De$ is an orbifold; 
\item $K_{\ov{X}_\De}= \pi^{*} K_{X_\De}$, i.e. the resolution is crepant. 
The proper transform $\ov{Z}_\De$ of $Z_\De$ is a smooth Calabi-Yau threefold and is a section of the base point free linear system $\left|-K_{\ov{X}_\De}\right|$.
\end{itemize}
\ePro
\noindent
Throughout this paper, by a Batyrev Calabi-Yau threefold, we mean a smooth generic Calabi-Yau threefold $\ov{Z}_\De$ given by Proposition~\ref{pro:refinements}.

There is a explicit formula for the hodge numbers of $\ov{Z}_\De$; see \cite{BC}. 
Since we are only interested in the K\"{a}hler cone of $\ov{Z}_\De$, we will only recall the calculation of $h^{1,1}(\ov{Z}_\De)$.

Let $\Theta$ be an arbitrary face of $\De$. Corresponding to $\Theta$, there is a dual face $\Theta^*$ of the dual polytop $\De^*$ defined by, 
$$ \Theta^*= \left\{ y\in \De^* \subset N \mid \left\langle x,y\right\rangle = -1, ~\forall x\in \De \right\}$$
If $\dim~\Theta=i$, then $\dim~\Theta^*=3-i$. 
For every $\Theta$, let $\ell(\Theta)$ denote the number of integral points on $\Theta$ and $\ell^*(\Theta)$ denote the number of interior integral points on $\Theta$, i.e. those that do not lie on proper subfaces of $\Theta$. 

\bPro{pro:h11}\textnormal(\cite[Proposition 4.4.1]{B1}\textnormal)
Given a 4-dimensional reflexive polytop $\De$,
\bEqu{equ:h11}
 h^{1,1}(\ov{Z}_\De) = \ell(\De^*) - 5 - \sum_{\dim~\Theta^*=3} \ell^*(\Theta^*) + \sum_{\dim~\Theta^*=2} \ell^*( \Theta^*)\cdot \ell^*(\Theta)
 \eEqu
\ePro 
\noindent
We will describe a natural base for $H^{1,1}(\ov{Z}_{\De})$ which explains the right hand side of (\ref{equ:h11}). We will use this base in the proof of Theorem~\ref{thm:main}.\\

Excluding the origin, there are $k= \ell(\De^*) - 1$ integral points on $\De^*$, all of which lie on the boundary of $\De^*$. 
Lets name these points $\left\{v_1,\ldots,v_k\right\}$. 
All of these points appear as the vertices of a maximal triangulation of $\De^*$, and therefore correspond to a toric divisor of $\ov{X}_\De$.  
Note that $\Si'$ is  the fan over  the refinement of  $\De^*$, obtained via a MPCP desingularization. 
Lets denote the corresponding divisors by $\left\{D_1,\ldots, D_k\right\}$. 
If $v_i$ is a an interior point of some 3-dimensional face $\Theta^*$ of $\De^*$, then $D_i$ lies over a singular 0-skeleton of $X_\De$; therefore, it does not meet $\ov{Z}_\De$ at any point, because $Z_\De$ is disjoint from the singular 0-skeleton of $X_\De$. 
The number of such divisors is $\sum_{\dim \Theta^*=3} \ell^*(\Theta^*)$.
This count has to be be deducted from the set of divisors which contribute to $H^{1,1}(\ov{Z}_{\De})$. 
Assume $\left\{v_1,\cdots,v_l\right\}$, 
$$l=\ell(\De^*) - 1 - \sum_{\dim \Theta^*=3} \ell^*(\Theta^*),$$ 
are all integral points of boundary of $\De^*$, except those which lie in the interior of a 3-dimensional face, and let $\left\{D_1\cdots D_l\right\}$  be the corresponding divisors. 
Every such divisor intersects $\ov{Z}_\De$ in a non-empty divisor; moreover, the intersection divisor is connected if and only if the corresponding point $v_i$ does not lie in the interior of a 2-dimensional face. 
In fact, if $v_i$ is an interior point of a 2-dimensional face $\Theta^*$, then $\Theta$ is 1-dimensional and $Z_\De \cap V_{\Theta}$ is made of $d(\Theta)= \ell^*(\Theta)+1$ singular points; therefore, $D_i \cap \ov{Z}_\De$ has $\ell^*(\Theta)+1$ irreducible components. 
This explains the factor $\sum_{\dim \Theta^*=2} \ell^*( \Theta^*)\cdot \ell^*(\Theta)$ in the right hand side of (\ref{equ:h11}). 
Considering all the irreducible components of the intersection of $\ov{Z}_\De$ and the toric divisors of $\ov{X}_\De$, we obtain 
\bEqu{equ:preset}
\ell(\De^*) - 1 - \sum_{\dim \Theta^*=3} \ell^*(\Theta^*) + \sum_{\dim \Theta^*=2} \ell^*( \Theta^*)\cdot \ell^*(\Theta)
\eEqu
irreducible divisors in $\ov{Z}_\De$, which has four divisors more than the right hand side of (\ref{equ:h11}). 
Every point $m\in M$ corresponds to rational function on $\ov{X}_\De$, and by restriction, to a rational function on $\ov{Z}_\De$. 
This way we get a four-dimensional space of linear relations between divisors  which has to be subtracted from (\ref{equ:preset}) to get the number of linearly independent divisors.
Summarizing, we get: 
\begin{itemize}
	\item Let $\left\{v_1,\ldots, v_a\right\}$ be the integral points on the boundary of $\De^*$ which do not lie on the interior of any two-dimensional or 3-dimensional face. 
Corresponding to these points, there are $a$ irreducible divisors $E_1,\ldots,E_a$, $E_i=D_i\cap \ov{Z}_\De$ , on $\ov{Z}_\De$.
	\item Let $\left\{u_1\cdots u_b\right\}$ be the integral points on the boundary of $\De^*$ which lie on the interior of some two-dimensional face $\Theta$. 
Corresponding to $u_i$, $1\leq i \leq b$, there are $d=d(\Theta)$ irreducible divisors on $\overline{Z}_{\De}$, $F_{i1},\ldots, F_{id}$, such that $\sum_{j} F_{ij}=D_i\cap \ov{Z}_\De$.
	\item The set of irreducible divisors $E_i$ and $F_{ij}$ generate $H^{1,1}(\ov{Z}_\De,\Q)$, and there is a four-dimensional set of linear relations among them; eliminating some four of these divisors, we get a basis of $H^{1,1}(\ov{Z}_{\De},\Q)$. 
\end{itemize}
In the next section, we will use these notations in the proof of theorem \ref{thm:main}.

\section{Second Chern class of Batyrev Calabi-Yau threefolds}\label{sec:main}
We start this section with an inspiring theorem of Kollar, which applies to $\ov{Z}_\De=Z_\De$ if $X_\De$ is smooth.

\bThm{thm:Kollar}\textnormal(\cite[appendix]{BK}\textnormal)
Given $X$, a smooth Fano variety with  $\dim_\C X\geq 4$, let $ Z\subset X$ be a smooth Calabi-Yau divisor in $\left|-K_X\right|$. Then, the natural map 
$\iota^*\colon \ov{\ck}_X \to \ov{\ck}_Z$ (or dually the map $\iota_*\colon \ov{NE}(Z) \to \ov{NE}(X)$) is an isomorphism of K\"{a}hler cones (cone of effective curves, respectively). 
Moreover, K\"{a}hler cone of $X$ and thus $Z$ is rational polyhedral.
\end{theorem}

The starting point of the proof is the Lefschetz Hyperplane Theorem, $H^{1,1}(X)\cong H^{1,1}(Z)$. 
For general $\De$,  $X_{\De}$ is a Fano 4-fold, but not all of them are smooth; therefore, a priori, we do not have a comparison between $H^{1,1}(X_\De)$ and $H^{1,1}(Z_\De)$. 
If  we move to a MPCP resolution $\ov{X}_{\De}$, 
Batyrev Calabi-Yau 3-folds $\ov{Z}_\De$ are still sections of $\left|-K_{\ov{X}_{\De}}\right|$, 
but $\ov{X}_\De$ is no longer Fano. 
Moreover, some divisors on $\ov{X}_\De$ do not intersect $\ov{Z}_\De$, and some others intersect in more than one component.
Followings is a example where this phenomena happens.

\begin{example}
Consider the fan $\Si$, given by following $15$ top dimensional cones in $\Z^4$, where $e_i$'s are a base of $N=\Z^4$:
\begin{enumerate}
\item $\si=\left\langle e_1,e_2,e_3,e_4,e_1+e_2+e_3-2e_4\right\rangle$;
\item $\left\langle \pm e_1, \pm e_2 , \pm e_3, e_4\right\rangle$ where all signs are not $+$ at the same time;
\item $\left\langle \pm e_1, \pm e_2, \pm e_3, e_1+e_2+e_3-2e_4 \right\rangle$, where all signs are not $+$ simultaneously.
\end{enumerate}
This is Gorenstein Fano fan and the corresponding polytop $\De$ is given by convex hull of following integral points in $M\cong \Z^4$:
\begin{enumerate}
\item $(\ep_1, \ep_2, \ep_3, -1)$ where $\ep_i= \mp 1$; $(-1,-1,-1,-1)$ corresponds to the first item above and rest of them correspond to the second group of the cones.
\item $-(\ep_1, \ep_2, \ep_3, \frac{\sum \ep_i -1}{2})$; where $\epsilon_i=\pm 1$ and all of them are not $+1$ simultaneously. These correspond to the 3rd group of the cones. 
\end{enumerate}
Let $X_\De$ be the corresponding toric variety. 
Among all the cones listed above, the first one is singular and its MPCP resolution $\ov{X}_\si=X_{\si'}$ is obtained by adding a wall 
(a small resolution, producing a rational curve $C$, above the corresponding singular point). 
After this small resolution, it still has a simplicial singular cone isomorphic to a  quotient singularity $\C^4/\Z_2$; thus, does not have any MPCP desingularization. 
All 3-dimensional subcones of the first group are non-singular; therefore, there is only one singular point on corresponding affine chart. 
All the cones in second the second group are smooth. 
The cones in 3rd group have quotient singularity $\C^2/\Z_2$; thus, again, they do not have any MPCP desingularization. 
We conclude that there is a total of 8 singular points on $\ov{X}_{\De}$, all of them isomorphic to $\C^4/\Z_2$. 
A generic Calabi-Yau 3-fold $Z_\De \subset X_\De$ is smooth and does not change through the MPCP-resolution $\ov{X}_\De \to X_{\De}$; i.e. $\ov{Z}_\De\cong Z_\De$.

Let $v=e_1+e_2+e_3-2e_4$, then any Weil divisor on $X_\De$ is of the form 
$$\sum_{i=1}^{3} ( a_i D_{e_i}+ b_i D_{(-e_i)})  + c D_{e_4}+ d D_v;$$ 
but some of them are not $\Q$-Cartier. 
The  $\Q$-Cartier divisors are of the form 
$$\sum_{i=1}^{3} (a_i D_{e_i}+ b_i D_{(-e_i)}) + c D_{e_4}+ (\sum a_i - 2c) D_v;$$ 
therefore, for $\Q$-Cartier divisors, the parameter $d$ can not be chosen arbitrarily. 
We conclude that the rational Picard group $Pic(X_{\De})\otimes \Q $ is isomorphic to $\Q^3$. 
Over $\ov{X}_\De$, all coefficients can be determined independently, i.e. $Pic(\ov{X}_{\De})= \Q^4$.
In fact, the proper transform of any Weil divisor to $\ov{X}_\De$ is $\Q$-Cartier.
From (\ref{equ:h11}) we get $h^{1,1}(Z_{\De})= 4 $. 
Considering the cone of effective curves, we observe that the rational curve $C$ obtained from resolution of $X_\De$ is in $\ov{NE}(\ov{X}_\De)$, but it can not be in the image of $\iota_* \ov{NE}(Z_\De)$, contrary to the smooth case of Theorem~\ref{thm:Kollar}.
This shows that $Pic(\ov{X}_\De)\otimes \Q \to Pic(Z_\De)\otimes \Q$ has non-trivial kernel and cokernel,
and $Pic(X_\De)\otimes \Q \to Pic(Z_\De)\otimes \Q$ is a proper embedding; see \cite[Proposition 3.5]{DKH}. 
\end{example}

For every smooth Calabi-Yau threefold $Z$, the second Chern class of the tangent bundle, $c_2(TZ)$, is an integral class in $H^4(Z)$, or dually in $N_1(Z)$; therefore, it corresponds to a linear map 
$c_2(TZ)\wedge - : H^{1,1}(Z)\to \Z$. By \cite[Theorem 1.1]{Mi}, we know that this linear map is non-negative on the K\"{a}hler cone, i.e. $\tn{PD}[c_2(TZ)] \in \ov{NE}(Z)$ .
If  $\tn{PD}[c_2(TZ)] \in NE^0(Z)$, where $NE^0(Z)$ is the interior of cone of effective curves, we get an intrinsically determined projective embedding of $Z$: since the intersection of $\ov{\ck}_{Z}$ with level sets $c_2(TZ)^{-1}(a)$, $a\in \Z$, are bounded, there are finitely many integral points on each level set. 
Let $a$ be the first positive integer where this set is non-empty. 
Consider the ample class obtained by sum of elements of this set and the projective embedding determined by a multiple of this class into some $\P^N$. 
Every automorphism of $Z$ comes from a linear automorphism of $\P^N$; therefore, $\tn{Auto}(Z)$ has finitely many components. Since it is discrete, it should be finite. 

\begin{proof}(of theorem \ref{thm:main})
Let $\De$ be a four-dimensional reflexive polytop, and $\ov{Z}_\De$ be a generic smooth Calabi-Yau 3-fold obtained from a MPCP resolution of $X_\De$. 
Consider the set of divisors $E_i$ and $F_{ij}$ listed at the end of Section~\ref{sec:basic-facts}.
If  $\ov{X}_{\De}$ is smooth, then the Chern classes of tangent bundle are given by 
$$c(T\ov{X}_{\De})=\prod(1+D_i),$$
where the product is over all of the boundary divisors. 
Since $\ov{Z}_\De$ is a smooth section of the base-point free linear system $\left|K_{\ov{X}_\De}\right|$, 
from adjunction formula we get,
\bEqu{equ:reduction}
\aligned
 c(T\ov{Z}_\De) &= \frac{c(T\ov{X}_\De)}{c(N_Z^{\ov{X}_\De})}= \frac{\prod (1+D_i)}{1+(\sum D_i)}\\
 &=\big(1+ (\sum D_i)+ (\sum_{i<j} D_i\cdot D_j)+\cdots)(1-(\sum D_i)+(\sum D_i)^2+\cdots\big)|_{\ov{Z}_\De}\\
 &=\big( 1+\;0+\; (\sum_{i<j} D_i\cdot D_j)- (\sum D_i)^2+ (\sum D_i)^2 + \cdots \big)|_{\ov{Z}_\De}.
 \endaligned
 \eEqu
Therefore,
$$c_2(T\ov{Z}_\De)=c_2(T\ov{X}_\De)= (\sum_{i<j} D_i\cdot D_j)\mid_{\ov{Z}_\De},$$ 
or equally
$$\tn{PD}[c_2(T \ov{Z}_{\De})]= \sum C_{ij},$$
where $C_{ij}= \ov{Z}_\De \cap D_i \cap D_j$ is a curve class on $\ov{Z}_\De$. 

More generally, if $\ov{X}_\De$ has isolated $\Q$-factorial singularities,  we use the following simple trick. Consider a toric resolution $Y\ra \ov{X}_\De$ whose exceptional divisors are supported on the isolated singularities of $\ov{X}_\De$. 
From the adjunction formula for $\ov{Z}_\De\subset Y$, we get
\bEqu{equ:reduction2}
 c(T\ov{Z}_\De) = 
 \frac{c(TY)}{c(N_Z^Y)}\cong 
 \frac{c(TY)}{c(N_Z^{\ov{X}_\De})}= 
 \frac{\prod (1+D_i^Y)}{1+(\sum D_i^{\ov{X}_\De})}, 
\eEqu  
where $D_i^Y$ are toric divisors on $Y$ and $D_i^{\ov{X}_\De}$ are toric divisors on $\ov{X}_\De$.
Since $\ov{Z}_\De$ does not meet the singularities of $\ov{X}_\De$, Equation~\ref{equ:reduction2} restricted to $\ov{Z}_\De$ is equal to Equation~\ref{equ:reduction}. Therefore, even if $\ov{X}_\De$ is not smooth,  Equation~\ref{equ:reduction2} is still true.
In order to finish to proof, 
we need to show that the set of curves $\{C_{ij}\}$ generate $H_2(\ov{Z}_\De,\Q)$.
Then every positive some of these effective curves lie in the interior of cone of effective curves. 

By Hard-Lefschetz theorem, we know that the map 
$$ H^{1,1}(\ov{Z}_\De,\Q) \otimes H^{1,1}(\ov{Z}_\De,\Q) \to H^4(\ov{Z}_\De,\Q) \cong H_2(\ov{Z}_\De,\Q),$$
is surjective. 
Therefore, elements of the form $E_i \cdot E_j$, $E_i \cdot F_{ij}$ and $F_{ij} \cdot F_{kl}$ generate $N_1(\ov{Z}_{\De})$. 
Right hand side of (\ref{equ:reduction2}) includes all these intersections,
except the intersection between two identical divisors. 
On the other hand, every toric divisor $D_i$ is equivalent to some linear combination of other toric divisors.
In fact, there is some $m\in M$ such that the linear relation given by $m$ express $D_i$ as a linear combination of several other toric divisors. 
Therefore, the intersection between identical divisors can be written as a linear combination of non-identical intersections. This shows that the set $\{C_{ij}\}$ of the effective curves in $\ref{equ:reduction2}$ includes a basis of  $H_2(\ov{Z}_\De,\Q)$.

\eProof

\bRem{rem:completeintersection}
Theorem~\ref{thm:main} is not true for complete intersections. 
The abelian fibered Calabi-Yau 3-fold example of \cite{GM} is the intersection of two hypersurfaces in $\P^2\times \P^2 \times \P^1$. Its Picard number is $19$ and the curve class corresponding to its second Chern class lies in a two-dimensional face of boundary. In this case, the autoharphism group is infinite, but has been verified to have a polyhedral fundamental domain on the K\"{a}hler cone.
\eRem


\end{document}